\title{Increasing the chromatic number of a random graph}
\author{Noga Alon 
\thanks{Schools of Mathematics and Computer Science, 
Sackler Faculty of Exact Sciences, Tel Aviv University, 
Tel Aviv 69978, Israel. Email: nogaa@tau.ac.il. 
Research supported in part by an ERC Advanced
grant,
by a USA-Israeli BSF grant, and by the Hermann Minkowski Minerva 
Center for Geometry
in Tel Aviv University.} \and
Benny Sudakov\thanks{Department of Mathematics,
UCLA,  Los Angeles, CA 90095. Email: {\tt bsudakov@math.ucla.edu}.
Research supported in part by NSF CAREER award DMS-0812005 and by a USA-Israeli BSF grant.}}

\documentclass[11pt]{article}
\usepackage{amsfonts,amssymb,amsmath,latexsym}

\newenvironment{proof}
      {\medskip\noindent{\bf Proof.}\hspace{1mm}}
      {\hfill$\Box$\medskip}

\def\qed{\ifvmode\mbox{ }\else\unskip\fi\hskip 1em plus 10fill$\Box$}

\newtheorem{theorem}{Theorem}[section]
\newtheorem{fact}[theorem]{Fact}
\newtheorem{lemma}[theorem]{Lemma}

\newtheorem{corollary}[theorem]{Corollary}

\newtheorem{prop}[theorem]{Proposition}

\newtheorem{definition}[theorem]{Definition}

\begin{document}
\date{}

\maketitle

\begin{abstract}
What is the minimum number of edges that have to be added to the
random graph $G=G_{n,0.5}$ in order to increase its chromatic number
$\chi=\chi(G)$ by one percent ? One possibility is to  add all
missing edges on a set of $1.01 \chi$ vertices, thus
creating a clique of  chromatic number $1.01 \chi$. 
This requires, with high probability,
 the addition of $\Omega(n^2/\log^2 n)$ edges.
We show that this is tight up to a constant factor, consider the 
question for more general random graphs $G_{n,p}$ with $p=p(n)$,
and study a local version of the question  as well.  

The question is motivated by the study of the resilience of graph
properties, initiated by the second author and Vu, and improves
one of their results.
\end{abstract}

\section{Introduction}
Consider the probability space whose points are graphs 
on $n$ labeled vertices, 
where each pair of vertices forms an edge,
randomly and independently with probability $p$. 
The random graph $G_{n,p}$ denotes a random point in this probability space.
This concept is one of the central notions in modern discrete
mathematics and it has been studied intensively during the last 50 years.
By now, there are thousands of papers and two excellent monographs 
by Bollob\'as \cite{Bol1} and by
Janson et al. \cite{JLR} devoted to 
random graphs and their diverse applications.
The subject of the theory of random graphs is
the investigation of the asymptotic behavior of various graph
parameters. We say that a graph property $\cal P$ 
holds {\em asymptotically almost surely}
(a.a.s.) if the probability that $G_{n,p}$ has $\cal P$ tends to one as
$n$ tends to infinity.

One of the most important parameters of the 
random graph $G_{n,p}$ is its
chromatic number, which we denote by 
$\chi(G_{n,p})$. Trivially for every graph
$\chi(G)\geq |V(G)|/\alpha(G)$, where $\alpha(G)$
denotes the size of the largest independent set in $G$. 
It can be easily shown, using first moment computations,
that a.a.s $\alpha(G_{n,p})\leq 2\log_b (np)$, where
$b=1/(1-p)$ (all logarithms in this paper are in the natural base $e$). 
This provides a lower bound on the chromatic number of the random graph, 
showing that
$\chi(G_{n,p})\geq \frac{n}{2\log_b (np)}$. 
The problem of determining the asymptotic behavior  of $\chi(G_{n,p})$, 
posed by Erd\H{o}s and R\'enyi in the early 60s, stayed for many years 
as one
of the major open problems in the theory of random graphs, until
its solution by Bollob\'as \cite{Bol2}, using a
novel approach based on martingales that enabled him to prove that
a.a.s. $\chi(G_{n,p})=(1+o(1)) \frac{n}{2\log_b (np)}$ 
for dense random graphs. Later \L uczak \cite{L} showed that this 
estimate also 
holds for all values of $p \geq c/n$. 
In this paper we strengthen these classical results, by showing that 
the chromatic number remains $(1+o(1)) \frac{n}{2\log_b (np)}$
even if we are allowed to add to $G_{n,p}$ any set of not 
too many additional edges.
To describe the main results it is convenient to use the framework of 
resilience, introduced by Sudakov and Vu \cite{SV}.

A graph property is called \textit{monotone
increasing (decreasing)} if it is preserved under edge addition
(deletion). Following \cite{SV}, we define:

\begin{definition}
Let $\mathcal{P}$ be a monotone increasing (decreasing) graph
property.
\begin{itemize}

\item 
The global resilience of $G$ with respect to
$\mathcal{P}$ is the minimum number $r$ such that by deleting (adding) $r$
edges from $G$ one can obtain a graph not having $\mathcal{P}$.

\item
The local resilience of a graph $G$ with respect
to $\mathcal{P}$ is the minimum number $r$ such that by deleting (adding) at
most $r$ edges at each vertex of $G$ one can obtain a graph not having
$\mathcal{P}$.

\end{itemize}
\end{definition}

Intuitively, the question of determining the resilience of a graph $G$
with respect to a graph property $\mathcal{P}$ is like asking, ``How
strongly does $G$ possess $\mathcal{P}$?''.
Using this terminology, one can restate many important 
results  in extremal graph theory in this language.
For example, the classical theorem of Dirac
asserts that the complete graph $K_n$ 
has local resilience $\lfloor n/2 \rfloor$ with respect to having 
a Hamilton cycle.
In \cite{SV}, the authors have
initiated the systematic study of global and local resilience of
random and pseudo-random graphs. 
They obtained resilience results
with respect to various properties such as perfect matching,
Hamiltonicity, chromatic number and having a nontrivial automorphism
(the last result appeared in an earlier paper with Kim \cite{KSV}). 
For example, they proved that if $p > \log^4 n /
n$ then a.a.s. any subgraph of $G(n,p)$ with minimum degree
$(1/2+o(1))np$ is Hamiltonian. Note
that this result can be viewed as a generalization of Dirac's
theorem mentioned above, since  a complete graph is also a random graph
$G(n,p)$ with $p=1$. This connection is natural
and most of the resilience results for random and pseudo-random graphs
can be viewed as generalizations of classical results from graph theory. 
For additional recent resilience type results, see, e.g.,
\cite{DKMS, FK, KLS, BKT,  BCS}.

In \cite{SV}, Sudakov and Vu proved that the local resilience of 
dense $G_{n,p}$ with respect to having chromatic 
number $(1+o(1))\frac{n}{2\log_b (np)}$ is at least $np^2/\log^5 n$.
The main aim of the present short  paper is to obtain the 
following new bounds on both the global and the local resilience of 
the chromatic number of the random graph, 
which substantially improve this result from \cite{SV}.

\begin{theorem}
\label{main1}
Let $\epsilon>0$ be a fixed constant and let $n^{-1/3+\delta} \leq p \leq 1/2$ for some $\delta>0$. Then a.a.s.
for every collection $E$ of $2^{-12}\epsilon^2\frac{n^2}{\log^2_b (np)}$ edges the chromatic number of $G_{n,p} \cup E$ is still 
at most $(1+\epsilon)\frac{n}{2\log_b (np)}$.
\end{theorem}
This shows that the global resilience of $G_{n,p}$ with respect
to having 
chromatic number at most \\
$(1+\epsilon)\frac{n}{2\log_b (np)}$ is of order
$n^2/\log^2_b (np)$. The result is tight up to a constant factor. Indeed, 
take an arbitrary set of, say, $n/\log_b (np)$ vertices of the 
random graph and add to it all the missing edges to make it a
clique. This adds only $\sim \frac{1}{4}n^2/\log^2_b (np)$
edges but increases the chromatic number by a factor of $2$.

\begin{theorem}
\label{main2}
Let $\epsilon>0$ be a fixed constant and let $n^{-1/3+\delta} \leq p \leq 1/2$ for some $\delta>0$. Then a.a.s. for every graph $H$ on $n$ vertices with maximum degree 
$\Delta(H) \leq 2^{-8}\epsilon \frac{n}{\log_b (np) \log \log n}$ the chromatic number of $G_{n,p} \cup H$ is still
at most $(1+\epsilon)\frac{n}{2\log_b (np)}$.
\end{theorem}
As before, by transforming a subset of $n/\log_b (np)$ vertices of
the random graph to
a clique, it follows that 
this result is tight up to the $\log \log n$ factor.
Both these theorems show that adding quite large and dense graphs to $G_{n,p}$ has very little impact on its chromatic number.
It may be instructive to compare the above two theorems 
to the following folklore result (see,
e.g., \cite{Lovaszbook} Chapter 9).

\begin{fact} 
Let $G$ and $H$ be two graphs on the same set of points. Then
$$\chi (G \cup H) \le \chi (G) \chi (H),$$
and there are pairs of $G$ and $H$ such that the equality  holds.
\end{fact}

The rest of this short paper is organized as follows. 
In the next section we prove our key technical lemma, which 
shows that in the random graph $G_{n,p}$ the independent 
sets of nearly maximal size are rather uniformly distributed. 
Using this lemma we establish Theorems \ref{main1} and \ref{main2} 
in Section 3.  The last section of the paper contains
some concluding remarks and open questions. Throughout 
the paper, we systematically omit floor and ceiling signs whenever 
they are not crucial for the sake of clarity of presentation.
We also do not make any serious attempt to optimize the absolute 
constants in our statements and proofs.

\section{The distribution of independent sets in random graphs}
In this section we prove the statement which will be  
our main technical tool.
It shows that no matter which set of $m$ edges we add to the 
random graph $G_{n,p}$, a.a.s. there will be an independent set of nearly
maximal size which contains only a few of these edges. In order to 
state our result precisely we need some preliminaries.

Let $ n^{-1/3+\delta} \leq p \leq 1/2$ for some $\delta>0$ and let $k_0=k_0(n,p)$ be defined by
\begin{equation}\label{k0}
k_0=\max\left\{k: {n\choose k}(1-p)^{{k\choose 2}}\geq
n^4\right\}\,.
\end{equation}
One can show easily that $k_0$ satisfies $k_0\sim 2\log_b(np)$
with $b=1/(1-p)$. Also, it follows from known results on the
asymptotic value of the independence number of $G(n,p)$ (see,
e.g., \cite{JLR}, \cite{AS}) that a.a.s. the difference between $k_0$
and the independence number of $G(n,p)$ is bounded by an absolute
constant, as long as $p(n)$ is in the above range.

Let $\mu$ be the 
expected number of independent sets of size $k_0$ in $G_{n,p}$.
Clearly  $$\mu ={n\choose k_0}(1-p)^{{{k_0}\choose 2}}\ge n^4,$$
by the definition of $k_0$. For a pair $u,v\in G_{n,p}$, 
let $Z_{u,v}$ be the random variable
counting the number of independent sets of size  $k_0$ in $G_{n,p}$ 
that contain both $u$ and
$v$. Let $\mu_0=E[Z_{u,v}]$, then 
$$ \mu_0={{n-2}\choose{k_0-2}}(1-p)^{{{k_0}\choose 2}}\ . $$ 
It is easy to see that $\mu_0/\mu=(1+o(1))k_0^2/n^2$.

Let $X$ be the random variable which is equal to the size of 
the largest collection of independent sets of size $k_0$ in the 
random graph $G_{n,p}$
such that no pair of vertices $u,v$ belongs to more than $4\mu_0$ of these sets. We need the following lemma, which shows that 
with high probability the value of $X$ is concentrated around $\mu$.

\begin{lemma}\label{le21}
$$Pr[X \le 3\mu/5]\le
e^{-\frac{\mu^2}{300\mu_0^2n^2p}}.$$
\end{lemma}

To prove this lemma we will need first to estimate from below the expectation of 
$X$. For a pair of vertices $u,v$ in $G_{n,p}$ set
$$ Z_{u,v}^+=\left\{\begin{array}{ll}
                    Z_{u,v}, & Z_{u,v}>4\mu_0 \geq 2\mu_0/(1-p),\\
                    0,       & \mbox{otherwise}\ .
                    \end{array}
              \right.
$$ 
We also define $Z^+=\sum_{u,v}Z_{u,v}^+$. This random variable has
been  considered before in \cite{KSVW}, 
where the authors studied the probability that the 
random graph $G_{n,p}$ contains an independent set of 
size $k_0$. We will need the following 
claim, proved in this paper. 

\begin{prop}\label{pr22}
$\mathbb{E}[Z^+]=o(\mu)$.
\end{prop}

From this proposition we can immediately deduce the following 
bound on the expectation of $X$.

\begin{corollary}
\label{cr23}
Let $X$ be the size of the largest collection of independent sets 
of size $k_0$ in the random graph $G_{n,p}$
such that no pair of vertices belongs to more than $4\mu_0$ of these sets. Then $\mathbb{E}[X]=(1-o(1))\mu$.
\end{corollary}

\proof
Let $\cal F$ be the collection of all independent sets in 
$G_{n,p}$ of size $k_0$. By the definition of $k_0$, we have that
$\mathbb{E}[|{\cal F}|]=\mu$. For every pair of vertices $u,v$ which is contained in more than
$4\mu_0$ independent sets of size $k_0$, delete all these sets from $\cal F$. Note that for every pair of vertices $u,v$ we deleted at most
$Z^+_{u,v}$ sets and therefore the total number of deleted sets is at 
most $Z^+$. It is easy to see that the remaining independent 
sets cover every pair of vertices at most $4\mu_0$  times and 
therefore their number is at most $X$.
By Proposition \ref{pr22}, this implies that 
$$\mathbb{E}[X] \geq\mathbb{E}[|{\cal F}|]-\mathbb{E}[Z^+]=(1-o(1))\mu.$$
This completes the proof, since clearly $\mathbb{E}[X] \leq \mu$. \qed

Let $\cal I$ be the largest collection of independent sets of size $k_0$ 
in the random graph $G_{n,p}$ such that no pair of vertices belongs 
to more than $4\mu_0$ of 
these sets and recall that $X=|{\cal I}|$. Note that when we
connect by a new edge a pair of non-adjacent vertices $u,v$ 
of $G_{n,p}$ we can decrease the value of $X$ 
only by the number of independent sets in $\cal I$ which contain $u,v$. By definition, this is at most $4\mu_0$. Now suppose we delete an existing 
edge $(u,v)$ of  the
random graph. Although this might create many new independent sets 
of size $k_0$, they  all contain $u,v$ and we can include only 
at most $4\mu_0$ of them in 
$\cal I$. Hence also in this case the value of $X$ changes by at 
most $4\mu_0$, i.e., $X$ is a so called $4\mu_0$-Lipschitz function. 
Now to finish the proof of Lemma 
\ref{le21} we apply a concentration inequality for such functions, 
proved by Alon, Kim and Spencer (\cite{AKS}, see also \cite{AS}, 
Theorem 7.4.3). They considered a random 
variable $Y$ given on the space generated by mutually independent $0/1$ choices such that probability that a choice $i$ is one is $p_i$. Let $c_i$ be such that 
changing choice $i$ can change $Y$ by at most $c_i$, $C=\max_i c_i$ 
and let the total variance satisfy 
$\sum_ip_i(1-p_i)c_i^2 \leq \sigma^2$. 
Then if $aC <2\sigma^2$ for some 
positive $a$, then 
$$Pr\big[Y-\mathbb{E}[Y]<-a\big] \leq e^{-a^2/(4\sigma^2)}.$$

\noindent
{\bf Proof of Lemma \ref{le21}.} As we already mentioned $X$ is a 
$4\mu_0$-Lipschitz random variable, 
which depends on ${n \choose 2}$ random choices for the edges 
of $G_{n,p}$. This implies that the total variance 
is at most 
$16\mu_0^2 {n \choose 2}p(1-p) \leq 8 \mu_0^2n^2p=\sigma^2$. 
Let $a =\mu/3$. 
Using that $k_0>1/p$ and $\mu/\mu_0=(1+o(1))n^2/k_0^2$, 
it is easy to check that $4\mu_0 a \leq 2\sigma^2=16  \mu_0^2n^2p$.
Note that by Corollary \ref{cr23}, if $X \leq 3\mu/5$ 
then we also have that $X-\mathbb{E}[X] <-a$. 
Therefore, the desired estimate for the lower tail of $X$ 
follows from the above inequality of \cite{AKS}.
\qed 

\noindent
{\bf Remark.}\, Using the same proof one can also obtain estimates on 
the upper tail of $X$ and in particular show that for any 
fixed $\delta>0$ 
$$Pr\big[\big|\,X-\mathbb{E}[X]\,\big| > \delta \mu \big] <2e^{-\frac{\delta^2\mu^2}{40\mu_0^2n^2p}}\,. $$

Finally we are ready to prove the main result of this section, 
which roughly says that in the random graph
$G_{n,p}$ independent sets of nearly maximal size are uniformly distributed in the following sense.
Suppose we add to the random graph some set $E$ of $m$ edges. Consider a random subset of vertices of size $k_0$. By averaging, we expect 
that $2m/n^2$ fraction of its pairs are edges from $E$. Our next lemma shows that with very high probability $G_{n,p}$ has an independent set of size $k_0$, 
which has this property up to a constant factor.

\begin{lemma}
\label{tool}
If $ n^{-1/3+\delta} \leq p \leq 1/2$ for some $\delta>0$, 
then with probability at least $1-e^{-n^{1+\delta}}$ the random graph 
$G_{n,p}$ has the following property. For every collection $E$ of $m$ edges, there is an independent set $I$ in $G_{n,p}$ of size $k_0$ such that
$I$ contains at most $7k_0^2\frac{m}{n^2}$ edges of $E$.
\end{lemma}

\proof
Let $\cal I$ be the largest collection of independent sets of size $k_0$ in random graph $G_{n,p}$ such that no pair of vertices belongs to more than $4\mu_0$ of
these sets. Recall that $\mu_0/\mu=(1+o(1))k_0^2/n^2$. Also, by 
the definition of $k_0$, 
it is easy to check that $k_0 \leq \frac{2}{p}\log n$. Together this implies that 
$$\frac{\mu^2}{300\mu_0^2n^2p} \geq \frac{n^2}{400k_0^4p} \geq \frac{n^2p^3}{10^4\log^4 n}>n^{1+\delta}.$$
Therefore, by Lemma \ref{le21} we have that with probability at least 
$1-e^{-n^{1+\delta}}$ the size of $\cal I$ is at least $3\mu/5$.
Consider an auxiliary bipartite graph $H$ with parts $\cal I$ and $E$ in which independent set $I \in {\cal I}$ is adjacent to edge $(u,v) \in E$ iff
both vertices $u,v$ belong to $I$. By the definition of $\cal I$, every edge $(u,v) \in E$ is contained in at most $4\mu_0$ sets from $\cal I$.
Therefore the number of edges $e(H)$ is bounded by $4\mu_0 m$. Thus there is an independent set $I \in {\cal I}$, whose degree in $H$ is at most
$e(H)/|{\cal I}|$. This $I$ contains at most 
$$\frac{e(H)}{|{\cal I}|}\leq \frac{4\mu_0 m}{|{\cal I}|}\leq \frac{20\mu_0 m}{3\mu}\leq 7k_0^2\frac{m}{n^2}$$ edges from $E$.
\qed

\section{Global and local resilience of the chromatic number}
In this section we prove our main results. We start by recalling 
several additional facts used in the proofs. 
The first is the following classical theorem of Tur\'an 
(see e.g., \cite{AS} p. 95), 
which provides a lower bound for the size of the maximum 
independent set in a graph.

\begin{lemma}
\label{turan}
Let $G$ be graph on $n$ vertices with $e(G)$ edges. Then the 
independence number $\alpha(G)$ of $G$  satisfies
$$\alpha(G) \geq \frac{n^2}{2e(G)+n}.$$
\end{lemma}

We also need the following simple lemma which estimates the number of edges spanned by small subsets of random graph.

\begin{lemma}
\label{density}
Let $n^{-1/3} \leq p \leq 1/2$ and let $\epsilon$ be a positive constant.
Then a.a.s. every subset of $G_{n,p}$ of size $ s \leq \frac{\epsilon n}{16\log (np)}$ contains
at most $\frac{\epsilon np}{8\log (np)}s$ edges.
\end{lemma}

\proof
Define $r=\frac{\epsilon n}{16\log (np)}$. The probability of existence of a subset violating the assertion of the lemma is at most
$$
\sum_{s=rp}^r{n \choose s}{{s \choose 2} \choose 2srp}p^{2srp} \leq 
\sum_{s=rp}^r \left(\frac{en}{s} \Big(\frac{es}{4rp}\Big)^{2rp}p^{2rp}\right)^s \leq 
\sum_{s=rp}^r \left(n \Big(\frac{e}{4}\Big)^{2rp}\right)^s=o(1).
$$
Here we used that $rp>\sqrt{n}$ together with the well known fact that ${a \choose b} \leq (ea/b)^b$.
\qed

\vspace{0.15cm}

Finally, recall the simple fact that any graph with chromatic 
number at least $r$ must have a subgraph with minimum degree $r-1$.

\vspace{0.15cm}
\noindent
{\bf Proof of Theorem \ref{main1}:}\, 
Let $E$ be an arbitrary set of $m=2^{-12}\epsilon^2n^2/\log^2_b (np)$ 
edges, suppose it has been added to $G_{n,p}$, and put $G=G_{n,p} \cup E$.
Let $b=1/(1-p)$ and let $k=2\log_b (np/\log^3n)=(1+o(1))2\log_b (np)$.
Since every subset of $G_{n,p}$ is a random graph itself and the number of subsets is at most $2^n$,
using Lemma \ref{tool} together with the union bound, we obtain that a.a.s. every subset of $G_{n,p}$ of size 
$s \geq n/\log^2 n$ has an independent set $I$ of size at least $k$ which contains at most
$7k_0^2\frac{m}{s^2} \leq 8k^2\frac{m}{s^2}$ edges of $E$. Repeatedly apply the following procedure until the remaining graph has at most
$\frac{\epsilon n}{16\log (np)}$ vertices. Given current subset, which has $s$ vertices,
find in it an independent set $I$ of $G_{n,p}$ of size at least $k$ which contains at most
$8k^2\frac{m}{s^2}$ edges of $E$. Apply Lemma \ref{turan} to the induced subgraph $G[I]$ to find an independent set of
$G$ of size $\frac{k^2}{16k^2m/s^2+k}$. Color it by a new color, remove its vertices from $G$ and continue. 

If the current subgraph of $G$ has  $s \geq 2^{-i}n$ vertices, then we can find in it 
an independent set of size  at least $\frac{k^2}{2^{2i+4}k^2m/n^2+k}$. Therefore if we start with at most $2^{-i+1}n$ vertices, then after pulling out 
$$\chi_i \leq 2^{-i}n\Big/\frac{k^2}{2^{2i+4}k^2m/n^2+k}=2^{i+4}\frac{m}{n}+2^{-i}\frac{n}{k}$$
independent sets we will remain with less than $2^{-i}n$ vertices. 
Let $i_0$ be such that $2^{i_0}=16\epsilon^{-1}\log (np)$.
Summing up for all $1 \leq i \leq i_0$ we obtain that we can color all but $\frac{\epsilon n}{16\log (np)}$ vertices of $G$ using only 

\begin{eqnarray*}
\sum_{i=1}^{i_0} \chi_i&=&\sum_{i=1}^{i_0} \left(2^{i+4}\frac{m}{n}+2^{-i}\frac{n}{k}\right) \leq 2^{i_0+5}\frac{m}{n}+n/k
\leq 2^9 \epsilon^{-1}\log (np) \frac{2^{-12}\epsilon^2n}{\log^2_b (np)}+n/k\\
&\leq &\frac{\epsilon}{4}\cdot \frac{\log (np)}{\log_b (np)}\cdot n/k+n/k \leq  
(1+\epsilon/4)n/k
\end{eqnarray*}
colors. 
Next we prove that the remaining $\frac{\epsilon n}{16\log (np)}$ vertices of $G$ can be colored by at most $r=\frac{\epsilon n}{3\log_b (np)}$ additional colors.

Indeed suppose that the remaining vertices form a graph with chromatic number more than $r$. Then this graph contains a subgraph $G'$ with 
minimum degree at least $r-1$. Denote by $s$ the number of vertices of $G'$ and note that $ r \leq s \leq \frac{\epsilon n}{16\log (np)}$ vertices. The number of edges in $G'$ is at least 
$s(r-1)/2$. On the other hand, using that $\frac{1}{p}\log (np) \geq \log_b (np)$ together with Lemma \ref{density}, we conclude that the number of edges in 
$G'$ is at most 
$$\frac{\epsilon np}{8\log (np)}s+|E| \leq \frac{\epsilon n}{8\log_b (np)}s+ \frac{2^{-12}\epsilon^2n^2}{\log^2_b (np)}\leq
\frac{\epsilon n}{8\log_b (np)}s+\frac{2^{-10}\epsilon n}
{\log_b (np)}s < s \frac{\epsilon n}{7\log_b (np)}< \frac{s(r-1)}{2}.$$
This contradiction shows that any $\frac{\epsilon n}{16\log (np)}$ vertices of $G$ can be colored by at most $\frac{\epsilon n}{3\log_b (np)}$colors.
Therefore the  chromatic number of $G$ is at most
$(1+\epsilon/4)n/k+\epsilon n/(3\log_b n)<(1+\epsilon)\frac{n}
{2\log_b (np)}$, completing the proof.
\qed

\vspace{0.15cm}
\noindent
{\bf Proof of Theorem \ref{main2}:}\,
The proof of this result follows along the same lines as the one for 
global resilience.
Let $H$ be a graph with maximum degree $\Delta \leq 2^{-8}\epsilon\frac{n}{\log_b (np) \log \log n}$ and let $G=G_{n,p} \cup H$.
Set $b=1/(1-p)$ and  $k=2\log_b (np/\log^3n)=(1+o(1))2\log_b (np)$. Let $S$ be a subset of $G_{n,p}$ of size
$s=|S| \geq n/\log^2 n$ and let $e(H[S])$ be the number of edges of $H$ spanned by $S$. Since $H$ has bounded maximum degree we 
have that $e(H[S]) \leq \Delta s/2$. We can again assume, by Lemma \ref{tool}, that
every such subset $S$ has an independent set $I$ of size at least $k$ which contains at most
$$7k_0^2\frac{e(H[S])}{s^2} \leq 8k^2\frac{e(H[S])}{s^2} \leq \frac{4k^2 \Delta}{s}$$
edges of $H$. Applying Lemma \ref{turan} to the induced subgraph $G[I]$, we find in it 
an independent set of $G$ of size $\frac{k^2}{8k^2\Delta/s+k}$.
Repeatedly color every such independent set by a new color and remove it from $G$ until the remaining graph has at most
$\frac{n}{\log^2 n}$ vertices. 

If the current subgraph of $G$ has  $s \geq 2^{-i}n$ vertices, then we can find in it
an independent set of size  at least $\frac{k^2}{2^{i+3}k^2\Delta/n+k}$. Therefore if we start with at most $2^{-i+1}n$ vertices, then after pulling out
$$\chi_i \leq 2^{-i}n\Big/\frac{k^2}{2^{i+3}k^2\Delta/n+k}=8\Delta+2^{-i}\frac{n}{k}$$
independent sets we will remain with less than $2^{-i}n$ vertices. 
Let $i_0$ be such that $2^{i_0}=\log^2 n$, then $i_0 \leq 3\log \log n$.
Summing up for all $1 \leq i \leq i_0$ we obtain that we can color all but $\frac{n}{\log^2 n}$ vertices of $G$ using only
$$\sum_{i=1}^{i_0} \chi_i = \sum_{i=1}^{i_0} \left(8\Delta+2^{-i}\frac{n}{k}\right) \leq \frac{n}{k}+8\Delta i_0 
\leq \frac{n}{k}+24 \log \log n \frac{2^{-8}\epsilon n}{ \log_b (np) \log \log n}\leq (1+\epsilon/4)\frac{n}{k}$$
colors.

Consider any subset $S$ of $s \leq \frac{n}{\log^2 n}$ vertices of $G$. By Lemma \ref{density} it has 
at most $\frac{\epsilon np}{8\log (np)}s$ edges of $G_{n,p}$. The number of edges of $H$ inside $S$ is clearly at most
$\Delta s/2$. Therefore there are at most $\big(\frac{\epsilon np}{8\log (np)}+\Delta/2\big)s$ edges in the subgraph of $G$ induced by $S$.
Using that $\frac{1}{p}\log (np) \geq \log_b (np)$ we conclude that
$G[S]$ has a vertex of degree at most
$$2e(G[S])/s \leq \frac{\epsilon np}{4\log (np)}+\Delta < 
\frac{\epsilon n}{2\log_b (np)}+\frac{2^{-8} \epsilon n}{\log_b (np) \log \log n} 
< \frac{3\epsilon n}{5k}-1.$$
This shows that we can color the remaining $\frac{n}{\log^2 n}$ 
vertices in $3\epsilon n/(5k)$ colors and the whole graph $G$ by
$(1+\epsilon/4 +3\epsilon/5)\frac{n}{k}<(1+\epsilon)\frac{n}
{2\log_b (np)}$ colors.
\qed

\section{Concluding remarks and open problems}

We have studied the global and local resilience of random graphs
with respect to the property of having a chromatic number close to
its typical value.  Our bounds for global resilience are tight
up to a constant factor, but the ones for the local case are only
tight up to a $\log \log n$ factor. It seems plausible to
conjecture that the assertion of Theorem \ref{main2} holds even when the
$\log \log n$ term is omitted in the hypothesis. 
It is also
possible that Theorem \ref{main1} can be strengthened, and that the
most economical way to increase the chromatic number of the random
graph $G_{n,p}$  by a factor of $(1+\epsilon)$ is to  construct an
appropriately large  clique in it.  If this is the case, then one
has to add to $G_{n,p}$, a.a.s., more than
$n^2/(16 \log^2_b n)$ edges in order to increase its chromatic number by
a factor of $(1+\epsilon)$, for any fixed $\epsilon >0$ and
sufficiently large $n$. This remains open.  It may also be
interesting to estimate the minimum number  of edges that have to
be added to $G_{n,p}$ in order to
increase the  chromatic number  by a lower order term.
This is related to the question of estimating  the concentration of
the chromatic number of random graphs, and appears to be difficult.

\end{document}